\documentclass{amsart}
\usepackage{amscd,amssymb,graphicx}
\usepackage{tikz}
\usetikzlibrary{arrows}

\tikzstyle{Bvertex}=[circle, black, fill, draw, inner sep=0pt, minimum size=3pt]

\author[Fialowski]{Alice Fialowski}
\address{
Alice Fialowski\\
University of P\'ecs and\\
E\"otv\"os Lor\'and University\\
Budapest, Hungary} \email{fialowsk@cs.elte.hu, fialowsk@ttk.pte.hu}
\author[Penkava]{Michael Penkava $^\dagger$}
\address{
Michael Penkava\\
University of Wisconsin-Eau Claire\\
Eau Claire, WI 54702-4004} \email{penkavmr@uwec.edu}
\thanks{Research of the second author was funded by grants from the
University of Wisconsin-Eau Claire.\\
$^\dagger$ Michael Penkava passed away on September 6, 2020.}

\theoremstyle{definition}

\def \ph{\varphi}

\def \ra{\rightarrow}

\renewcommand{\hom}{\operatorname{Hom}}

\def \tns{\otimes}

\def \mcom{,\cdots,}
\def \k{\mbox{$\mathbb K$}}

\def \C{\mbox{$\mathbb C$}}
\def \R{\mbox{$\mathbb R$}}
\def \Z{\mbox{$\mathbb Z$}}

\def \sl{\mathfrak{sl}}
\def \so{\mathfrak{so}}

\def\sh{\operatorname{Sh}}

\def\im{\operatorname{Im}}

\def\linf{\mbox{$L_\infty$}}
\def\and{\mbox{ \rm and }}

\def\s#1{(-1)^{#1}}

\def\pha#1#2{\ph^{#1}_{#2}}

\def\psa#1#2{\psi^{#1}_{#2}}

\def\ho{\text{ho}}

\def\dinfty{\mbox{$d^\infty$}}

\input{xy}
\xyoption{all}

\newcommand{\oscl}{$\mathfrak g_{6,2}(\lambda)$}
\begin{document}
\setlength{\multlinegap}{0pt}
\title[Lie algebras with Invariant Inner Product]
{On the cohomology of Lie algebras with an invariant inner product}

\address{}%
\email{}%

\thanks{}%

\date{\today}
\begin{abstract}
In this work we consider all metric Lie algebras, having a nondegenerate symmetric invariant bilinear form, over $\C$ and $\R$ up to dimension 5 and all metric Lie algebras over $\C$ in dimension 6. We introduce cyclic and reduced cyclic cohomology to identify their metric deformations.
\end{abstract}
\maketitle

\section{Introduction}

Recently, Lie algebras with an invariant inner product became an intensive topic of research in Lie Theory.  Any reductive Lie algebra has such an invariant inner product, related to the Killing form, which is an invariant inner product on a semisimple Lie algebra.  Usually, one considers such algebras
over the real numbers, but complex forms also play a role. One advantage in considering the complex case is that in the complex case all inner products are equivalent, whereas over the real numbers the signature of the form also plays a role. Lie algebras with an invariant inner product are also referred to as \emph{metric Lie algebras} in the literature, and include the 4-dimensional diamond and oscillator algebras as examples of solvable real algebras with an invariant inner product.

Representations of the diamond Lie algebra have been studied \cite{cms,lpx,med-rev} but in general, the invariants of metric Lie algebras and the cohomology and deformations of such algebras has not been studied. 

There is a notion
of cyclic cohomology, which plays a role in the study of deformations preserving an invariant inner product.  We examine this notion in some depth in this paper, and apply a variant of this type of cyclic cohomology to study metric deformations of these metric algebras, that is deformations of these into other metric algebras. 

Associative algebras with invariant inner products and their cyclic cohomology were studied in \cite{ps2} and the connection between cyclic cohomology and deformations of these algebras preserving an invariant inner product was analyzed.  Later, in \cite{pen1,pen4}, a notion of cyclic cohomology of Lie algebras was studied and it was shown that cyclic cohomology classifies deformations preserving the invariant inner product.

An invariant inner product on a Lie algebra $V$ is a symmetric, nondegenerate bilinear form $\beta$ on $V$ which is invariant in the sense that 
\begin{equation*}
\beta([x,y],z)=\beta(x,[y,z])    
\end{equation*}
for all $x$, $y$ and $z$ where $[\,,\,]$ is the Lie bracket.  The existence of an invariant nondegenerate bilinear form determines a \emph{metric} Lie algebra. 

\section{Cohomology and Cyclic Cohomology}

Recall that the space of cochains $C(V,M)$ of a Lie algebra with coefficients in a module $M$ is given by $C(V,M)=\hom(\bigwedge,M)$. If $C^n(V,M)=\hom(\bigwedge^nV,M)$, then $C(V,M)=\Pi_{n=0}^\infty C^n(V,M)$. By convention, we usually consider $$C(V,M)=\bigoplus_{n=0}^\infty C^n(V,M).$$
(This difference is only important when considering \linf\ algebras, which are a generalization of Lie algebras.)

There is a map $D:C(V,M)\ra C(V,M)$, called the \emph{coboundary operator} which is determined by the maps $D:C^n(V,M)\ra C^{n+1}(V,M)$ given by
\begin{align*}
  D(\ph)(v_1\mcom v_{n+1})&=
  \sum_{\sigma\in\sh(1,n)}\s{\sigma}v_{\sigma(1)}.\ph(v_{\sigma(2)}\mcom v_{\sigma(n+1)})\\&+
  \sum_{\sigma\in\sh(2,n-1)}\s{\sigma}\ph([v_{\sigma(1)},v_{\sigma(2)}],v_{\sigma(3)},\mcom v_{\sigma_{n+1}}),
\end{align*}
where $\ph$ is an $n$-cochain, $\sh(k,\ell)$ is the set of \emph{(un)shuffles} of type $(k,\ell)$, that is those permutations of $k+\ell$ which are increasing on $1\cdots k$ and $k+1\cdots k+\ell$, $\s{\sigma}$ is the sign of the permutation $\sigma$ and  $[\cdot,\cdot]$ represents the Lie bracket on $V$.

A standard fact from Lie theory is that the map $D$  satisfies $D^2=0$ and so we can define the (\emph{Eilenberg-Chevalley}) cohomology
\begin{equation*}
  H^n(C,V)=\ker(D:C^n(V,M)\ra C^{n+1}(V,M))/\im(D:C^{n-1}(V,M)\ra C^n(V,M)).
\end{equation*}

The cochains $C(V,V)$ with coefficients in the \emph{adjoint representation} denote by $C(V)$, so that $C^n(V)=C^n(V,V)$ are the $n$-cochains of the Lie algebra, and $H^n(V)=H^n(V,V)$ is called the \emph{cohomology} of (the Lie algebra structure on) $V$.

An inner product $\langle\cdot,\cdot\rangle$ on $V$ is \emph{invariant} with respect to a Lie algebra structure if
\begin{equation*}
  \langle [u,v],w\rangle=\langle u,[v,w]\rangle.
\end{equation*}
Given any inner product, we can define a cochain $\ph\in C^n(V)$ to be \emph{cyclic} if the associated element $\tilde:\bigwedge V\tns V\ra\k$, given by
\begin{equation*}
  \tilde\ph(v_1\mcom v_{n+1})=\langle\ph(v_1\mcom v_n),v_{n+1}\rangle,
\end{equation*}
satisfies
\begin{equation*}
  \tilde\ph(v_{n+1},v_1\mcom v_n)=\s{n}\tilde\ph(v_1\mcom v_{n+1}).
\end{equation*}
Note that an inner product is invariant with respect to a Lie algebra structure $d\in C^2(V)$ precisely when $\tilde d$ is cyclic with respect to the inner product.  Also, $\ph$ is cyclic precisely when $\tilde\ph$ is antisymmetric, that is $\tilde\ph\in C^{n+1}(V,\k)$, where $\k$ is the trivial
module structure.

If $\ph$ and $\psi$ are cyclic cocycles, the Chevalley-Eilenberg bracket of the two cochains $[\ph,\psi]$ is also cyclic (see \cite{ps2,pen4}),
where if $\ph\in C^k(V)$  and $\psi\in C^l(V)$, then $[\ph,\psi]\in C^{k+l-1}(V)$ is given as follows.
Let \begin{align*}
(\ph\circ\psi)(v_1\cdots v_{k+l-1})=
\sum_{\sigma\in\sh(k-1,l)}(-1)^\sigma\ph(\psi(v_{\sigma(1)}\cdots v_\sigma(l))v_{\sigma(l+1)}\cdots v_{\sigma(k+l-1)}),
\end{align*}
then \begin{equation*}
[\ph,\psi]=\ph\circ\psi-(-1)^{(k+1)(l+1)}\psi\circ\ph.  
\end{equation*}
This bracket equips $C(V)$ with the structure of a $\Z$-graded superalgebra. In particular, the $\Z$-graded super Jacobi identity holds, as well as the graded antisymmetry, given by 
\begin{equation*}
[\psi,\ph]=\s{(k+1)(l+1)+1}[\ph,\psi].
\end{equation*}

This allows us to extend the Chevalley-Eilenberg bracket of cochains in $C(V)$ to $C(V,\k)$, by
defining
\begin{equation*}
  [\tilde\ph,\tilde\psi]=\widetilde{[\ph,\psi]},
\end{equation*}
which equips $C^k(V,\k)$ with a graded Lie algebra structure. In fact, we can compute that if $\ph\in C^k(V)$, $\psi\in C^\ell(V)$, then
\begin{equation*}
  [\tilde\ph,\tilde\psi]=\sum_{\sigma\in\sh(l,k)}\tilde\ph(\tilde\psi(v_{\sigma(1)}\mcom v_{\sigma(\ell)}),v_{\sigma(\ell+1)}\mcom v_{\sigma(k+\ell)}).
\end{equation*}
It is not obvious why this formula is correct, since it seems unbalanced in terms of $\ph$ and $\psi$, but it can be checked that
with this bracket
\begin{equation*}
  [\tilde\psi,\tilde\ph]=\s{(k+1)(\ell+1)+1}[\tilde\ph,\tilde\psi],
\end{equation*}
which is the graded antisymmetry corresponding to the antisymmetry of the bracket of $\ph$ and $\psi$.
Note that the coboundary operator on an element $\tilde\ph\in C^n(V,\k)$ is just given by
\begin{align*}
  (D\ph)(v_1\mcom v_{n+2})&=\sum_{\sigma\in\sh(2,n)}\s{\sigma}\tilde\ph([v_{\sigma(1)},v_{\sigma(2)}],v_{\sigma(3)}\mcom v_{\sigma(n+2)})\\&=[\tilde\ph,\tilde d](v_1\mcom v_{n+2})
\end{align*}
This means that if we define the coboundary operator $D$ on the space $CC^n(V)$ of cyclic cochains by $$D(\tilde\ph)=[\tilde\ph,\tilde d],$$ then the associated cohomology
$HC(V)$, called the \emph{cyclic cohomology} of $V$, coincides with the shifted cohomology $H(V,\k)$ of $V$ with coefficients in $\k$.  That is, we have
$HC^n(V)=H^{n+1}(V,\k)$, for $k>0$. When $k=0$, there are no $-1$-degree cyclic cochains, but there can be $0$-degree cochains in $C(V,\k)$, so $HC^0(V)$ may not coincide with $H^1(V,\k)$. In fact, $HC^0(V)=Z^1(V,\k)$, because we ignore the 0-coboundaries in computing the cohomology.

When there is no invariant inner product, the cohomology $H(V,\k)$ is still well defined, and it is natural to define the cyclic cohomology of $V$ by
$$HC^n(V)=\begin{cases}
Z^1(V,\k),& n=0\\
H^{n+1}(V,\k),& n\ge 1
\end{cases}
$$ so that the definition of cyclic cohomology can be given independently of any invariant inner product.  
This is similar to the
associative algebra case, where one can define cyclic cochains without reference to an invariant inner product. But it should be noted that the bracket of
cyclic cochains makes sense only in the presence of an invariant inner product.

In the presence of an invariant inner product, there is a connection between cyclic cohomology and deformations of an algebra preserving the invariant inner product. However, as we shall see on a simple example, the connection is not as straightforward as one might think.

We introduce a notation for an element of $C^n(V)$ by defining $\ph^I_i(e_J)=\delta^I_J e_i$, where $I,J$ are strictly increasing multi-indices
of length $n$; i.e.  $I=(i_1\mcom i_n)$ with $i_k<i_{k+1}$ for $k=1\cdots n-1$. This allows us to represent an algebra as a 2-cochain in a very simple
form. By convention, we use the symbol $\psi$ in place of $\ph$ to  represent multi-indices which are of even length, and reserve the symbol $\ph$ for multiindices of odd length, corresponding to the fact that when $I$ is of even  length, the associated map $\psi^I_i$ is an \emph{odd} cochain and when $I$ is of odd length, the associated cochain is an even element in the superalgebra structure on the space of cochains. 

In particular, a 2-cochain can always be expressed in the form 
\begin{equation*}
    \psi=\psi^{i,j}_kc^k_{i,j},
\end{equation*}
using the Einstein summation convention, where $i<j$ and $c^k_{i,j}$ are the structure constants determining the Lie algebra. In terms of the standard form for describing the Lie bracket, the above formula translates to 
\begin{equation*}
    [e_i,e_j]=c^k_{i,j}e_k.
\end{equation*}
Moreover, the Jacobi identity for an algebra given by the cochain $\psi$ translates into the codifferential equation
\begin{equation*}
    [\psi,\psi]=0.
\end{equation*}

Consider the simple three dimensional complex Lie algebra $\sl(2,\C)$ on $V=\langle e_1,e_2,e_3\rangle$, which can be given by the cochain $d=\psa{12}3+\psa{23}1+\psa{13}2$ (in the classical notation, the nontrivial brackets are $[e_1,e_2]=e_3, [e_2,e_3]=e_1, [e_1,e_3]=e_2$). It has an invariant inner product represented by the
identity matrix. It is well known that $H^1(V,\C)$ and $H^2(V,\C)$ both vanish and that $H^3(V,\C)$ is 1-dimensional.  But this means that
$HC^2(V)=\langle d\rangle$ is generated by the cochain represented by the algebra itself, and that means that this cochain is \emph{nontrivial}. 
It is not difficult to show that  $d$ is not the coboundary of a \emph{cyclic} 1-cochain, which means that it really is true that the cyclic cohomology $HC^2(V)$ has an extra basis
element. This means that in some sense, the algebra deforms along itself, which would not make any sense in the usual notion of deformation of an algebra.  

In fact, we will explain later that the two algebras $d$ and $(1+t)d$, while isomorphic, are not formally isomorphic in the metric sense, and that is the source of the problem, because deformation theory considers a deformation to be trivial only when it is generated by a formal isomorphism.  The problem turns out to be related to the fact that the identity matrix is never a cyclic 1-cochain, and for reductive algebras, a multiple of the identity matrix determines the only trivial deformation taking $d$ to $(1+t)d$.

The problem with the cohomology can be analyzed as follows:
\begin{align*}
  \ker(D:CC^n\ra CC^{n+1})&=\ker(D:C^n\ra C^{n+1})\bigcap CC^{n}\\
  \im(D:CC^{n-1}\ra CC^n)&\subseteq\im(D:C^{n-1}\ra C^n)\bigcap CC^{n}
\end{align*}
The inclusion on the bottom may be strict.  One idea is to replace the left hand side with the right hand side, and define
a new type of cohomology, which we call \emph{reduced cyclic cohomology}
$$HRC^n=\ker(D:C^n\ra C^{n+1})\bigcap CC^{n}/(\im(D:C^{n-1}\ra C^n)\bigcap CC^{n}).$$ If we do this, then
we obtain that $HRC^2(V)$ vanishes, since $d$ is a coboundary of an ordinary 1-cochain. 
In fact, it is always true for any Lie algebra $d$ that $d=[d,I]$, where $I\in\hom(V,V)$ is the identity map, but the identity map is never a cyclic 1-cochain.  Note that the invariant inner product was used in the definition of reduced cyclic cohomology. In fact, if the algebra is reductive (or more generally, contains a simple quotient), then $HC^2$ does not vanish. However, for the solvable algebras that we have studied, we have determined that the reduced and ordinary cyclic cohomology coincide, which is a function of the type of invariant inner products that arise in those cases.


The reduced cyclic cohomology expresses deformations of metric Lie algebras if we allow formal deformations to be induced by any 1-cochain, instead of restricting to only cyclic cochains. This arises because we allow cyclic coboudaries of non cyclic cochains in the reduced cyclic cohomology, so we allow formal deformations which take a metric algebra to another metric formal algebra, rather than just metric formal deformations. However, there are problems with this point of view as well, because there is no reason to believe that the deformation given by the exponential of a non-cyclic 1-cochain $\beta$ preserves the invariant inner product, even when $[d,\beta]$ is a cyclic 2-cochain.

It is well known that the Killing form gives an invariant inner product for a semisimple Lie algebra, and thus reductive Lie algebras have an invariant
inner product.  However, these are not the only types of Lie algebras with invariant inner product.  Examples of real 4 dimensional solvable Lie algebras with
an invariant inner product are the diamond and oscillator algebras, while in dimension 5, the nilpotent real Lie algebra $W_3$, which we discuss later, also has
an invariant inner product.  These cases have been well studied \cite{ko1,ko2,bord,astr}.
We will study low dimensional examples of real and complex Lie algebras and their deformations preserving the invariant inner product.

\section{Deformations}

Recall that a 1-parameter deformation of a Lie algebra structure
$d$ is a formal power series of the form
$$d_t=d+t\psi_1+t^2\psi_2+\cdots,
$$
where $\psi_k\in C^k(V)=\hom(\bigwedge^k V,V)$ are \emph{2-cochains}, see \cite{fi2}. The connection with cohomology is given by the fact that $\psi_1$ is a 2-cocycle, and moreover, we have
$$
D(\psi_{n})=-\tfrac12\sum_{k+l=n}[\psi_k,\psi_l],
$$
where $[\psi_k,\psi_l]$ is the bracket of the cochains $\psi_k$ and $\psi_l$.
If $d_t$ is isomorphic to some algebra structure $d'$, then we say that $d_t$ and $d'$ are equivalent, and we write $d_t\sim d'$. 

If $d_t\sim d'$ for all $t$ in some neighborhood of the origin, then the deformation is called a \emph{jump deformation} from $d$ to $d'$. If, on the other hand, $d_t\not\sim d_{t'}$ for $t'\ne t$ in some neighborhood of the origin, then the deformation is called a \emph{smooth deformation}. In this case, the set of algebras $d_t$ form a family of nonisomorphic algebras.

Multiparameter deformations are also possible, and there is a special type of multiparameter deformation, called a \emph{versal deformation}, which is of the form 
$$d^\infty=d+t_i\delta^i+\text{higher order terms},
$$
where the expression $t_i\delta^i$ represents the Einstein summation notation for a sum of basis elements $\delta^i$ for the cohomology $H^2$.  There are relations, of the form 
$r_i=t_kt_lr^{kl}_i+\text{ho}$, which are formal power series of order at least 2, with the number of relations being equal  to the dimension of $H^3$. The base $A$ of the deformation is the \emph{formal} algebra $A=\k[[t_1\mcom]]/(r_i)$, that is, the quotient of the ring of formal power series over the field $\k$ by the ideal generated by the relations.

In fact, if we take $\langle\delta^i\rangle$ to be a basis of $H^2$, $\langle\gamma^i\rangle$ a basis of the 2-coboundaries, $\langle\alpha^i\rangle$ a basis of $H^3$, and $\tau^i$ a basis of the 4-coboundaries, then there is a unique choice of parameters $x_i$, given by power series in the parameters $t_i$ of order at least 2, such that the deformation $d^\infty$, given by 
$$
d^\infty=d+t_i\delta^i+x_i\gamma^i
$$
satisfies 
$$[d^\infty,d^\infty]=r_i\alpha^i+u_i\tau^i,$$
with both $r_i$ and $u_i$ being given by power series in the $t_i$ of order at least $2$.  Moreover, the $u_i$-s are contained in the ideal generated by the $r_i$-s, which are called the relations on the base of the deformation.

As a consequence, if we solve the relations on the base, we obtain solutions to the Lie algebra condition $[d^\infty,d^\infty]=0$.  
The deformation $d^\infty$ is called a 
\emph{miniversal} deformation of $d$, and all formal deformations of $d$ can be obtained from $d^\infty$ in a natural manner. 

In a series of articles on low dimensional algebras, the authors have computed versal deformations for all algebras up through dimension 5 over $\C$, and for dimension up to 4 over $\R$..

For cyclic cohomology, we obtain a similar picture of deformations, including a versal deformation, where this time, the $\delta^i$, $\gamma^i$, $\alpha^i$ and $\tau^i$ are cyclic cochains, with $\delta^i$ and $\alpha^i$ given by basis of $HC^2$ and $HC^3$ respectively.  In order to do the computations necessary for this paper, we had to adapt the routines we used to compute versal deformations of algebras to compute the versal cyclic deformations, and this turned out to be more involved.

The source of the problem with deformations arises because of the  notion of formal (or infinitesimal) equivalence of deformations.  Two infinitesimal deformations $d$ and $d'$ are formally equivalent if there is a linear map $\beta:g\ra g$ such that if $g=\exp(t\beta)$
then $$d'=g^*(d)=d+t[d,\beta].$$ 
Now, we want $d'$ to be a metric algebra with respect to the same inner product, so there are two ways to guarantee this.  We can require $\beta$ to be cyclic with respect to $d$, in which case $[d,\beta]$ is automatically cyclic, since the bracket of two cyclic cochains is cyclic, or we can just require that $[d,\beta]$ be cyclic. The latter case is potentially problematic for formal deformations, because in that case, higher order terms arise, which may not be cyclic.  Thus, the restriction that $\beta$ is cyclic is natural, and gives rise to a consistent deformation picture.  However, this is precisely what we obtain from the cyclic, rather than the reduced cyclic, cohomology.  This means that the algebras $d$ and $(1+t)d$ may not be formally equivalent, and this is exactly what happens for any simple algebra, because the identity map is never a cyclic cochain.  

If $d_1$ and $d_2$ are two formal deformations of a Lie algebra with respect to the usual cohomology, then they are said to be formally equivalent if there is a collection $\beta_i$ of $1$-cochains for $i=1,\dots$ such that the action of $\exp(t^i\beta_i)$ on $d_1$ yields $d_2$.  The question is what to do if we require from the cochains $\beta_i$ to be cyclic.  Clearly, if we require $\beta_i$ to be a cyclic cochain for all $i$, this will be sufficient to ensure that $\exp(t^i\beta_i)^*\psi$ is cyclic, whenever $\psi$ is a cyclic cochain. But is  the condition necessary?  The answer turns out to be no!  In fact,
if $\beta=I$ is the identity map and 
$$d_1=d+t^i\psi_i,$$
then 
$$\exp(tI)*d_1=d+t(\psi_1+d)+t^2(\psi_2+\psi_1)+\cdots t^k(\psi_k+\psi_{k-1})+\cdots,$$
which is cyclic because all of the $\psi$ terms are assumed to be cyclic.  
We will see how this applies in our examples below.

\section{Dimension 3}
The simple complex 3-dimensional Lie algebra is $\sl(2,\C)$. 
It can be given by $$d=\psa{12}3-2\psa{13}1+2\psa{23}2.$$ 
This is the only nontrivial 3-dimensional complex Lie algebra with an invariant inner product.  There are two real forms for this complex algebra, $\sl(2,\R)$ and $\mathfrak{so}(3,\R)$.

The first real form can be given by the structure $d$ above, with respect to some basis $\langle e_1,e_2,e_3\rangle$.
It has an invariant inner product with respect to this basis, given by the matrix $ \left[ \begin {array}{ccc} 0&1&0\\ \noalign{\medskip}1&0&0
\\ \noalign{\medskip}0&0&2\end {array} \right]$. The signature of this matrix is $(2,1)$, but the signature of a metric Lie algebra is only determined up to the transposition of the signature, so that there is an invariant inner product with signature $(1,2)$ as well.  However, we can give it uniquely by requiring the first number to be greater than or equal to the second.  In fact, the signature of the metric is not always determinate.  For example, for the 2-dimensional trivial algebra, any invertible matrix will serve as a metric, so the signature is not well defined by the algebra. 

The second real Lie algebra is $\mathfrak{so}(3,\R)$, which can be given by the
structure
$$d=\psa{12}3-\psa{13}2+\psa{23}1,$$ with invariant inner product given by the identity matrix, so its signature is $(3,0)$.

Note that it is not surprising that the invariant matrix is not the same as for $\sl(2,\R)$, because the form of the matrix depends on the choice of basis. However, here, since we are over $\R$, the signature of the matrix comes into play, so in fact, the two invariant inner products are not of the same type, since they have different signatures.

In the Table below, we give the dimensions of the cohomology for the simple 3-dimensional Lie algebras, for $HC^n$, $HRC^n$, and the standard cohomology $H^n$, with coefficients in the adjoint representation.

\begin{table}[ht]
  \begin{tabular}
    {lcccc}
    $n$&\vline&$HC^n$&$HRC^n$&$H^n$\\
    \hline
    0&\vline&0&0&0\\
    1&\vline&0&0&0\\
    2&\vline&1&0&0\\
    3&\vline&0&0&0\\
    \hline\\
  \end{tabular}
  \caption{Cohomology of the 3-dimensional simple algebra}
\end{table}

The Table applies for the complex simple Lie algebra and both of its real forms. The only variations will be for the basis elements of the cohomology.

While there are no deformations of $\sl(2,\C)$, or any of its real forms, because the cohomology $H^2$ vanishes, something interesting occurs with cyclic cohomology.  With respect to $HC^2$, we have a nontrivial  cocycle given by $d$ itself, and therefore, the 1-parameter deformation 
$d_t=(1+t)d$  is not a trivial deformation, which is very unusual, since this type of deformation would be trivial with respect to the usual notion of deformation. Note that this deformation preserves the metric.  However, it is easy to see that if we take $\beta=cI$, then
$$\exp(t\beta)^*(d)=\exp(ct)d,$$ and we can solve $\exp(ct)=1+t$ for $c$. Thus, the two deformations are equivalent if we allow the application of an exponential of a term which, although not cyclic, always produces a cyclic cochain.

\section{4-dimensional Lie algebras with invariant inner product}
Since cyclic cochains $CC^n(V)$ correspond to cochains in $C^{n+1}(V,\k)$, we have
$cc^n=\binom4{n+1}$, where $cc^n$ is the dimension of $CC^n$. Thus, for all 4-dimensional cyclic algebras, we have 
$$
cc^0=4,\qquad cc^1=6,\qquad cc^2=4\qquad cc^3=1,
$$
and these are the only degrees in which there are nonvanishing cochains.

\subsection{The direct sum $\sl(2,\C)\oplus \C$ and its real forms}
This algebra can be given by $d=\psa{12}3+2\psa{23}2-2\psa{13}1$, and in this form, an invariant inner product is given by the matrix
$$B= \left[ \begin {array}{cccc} 0&1&0&0\\ \noalign{\medskip}1&0&0&0
\\ \noalign{\medskip}0&0&2&0\\ \noalign{\medskip}0&0&0&1\end {array}
 \right].$$

We could calculate the cyclic cochains, coboundaries and cocycles by our computer technology, but here we think it is technically easier to use standard results on the computation of cohomology of Lie algebras.
First, we recall that 
$$H^n(\mathfrak g,M\oplus N)=H^n(\mathfrak g,M)\oplus H^n(\mathfrak g,N)$$ if $\mathfrak g$ is a Lie algebra and $M$ and $N$ are $\mathfrak g$-modules.
Secondly, we recall the K\"unneth formula:
$$
H^n(\mathfrak g\oplus \mathfrak h,M)=\bigoplus_{k+\ell=n}H^k(\mathfrak g,M)\tns H^\ell(\mathfrak h,M),
$$
where $M$ is both a $\mathfrak g$ and $\mathfrak h$-module. Finally, we use the isomorphism between cyclic cohomology and shifted ordinary cohomology with trivial coefficients (except for $n=0$).  

Finally, we use the fact that if $\mathfrak g$ is simple (or semisimple), then $H^n(\mathfrak g,\mathfrak g)=0$ for all $n$.
Let $\mathfrak g=s$ be $\sl(2,\C)$ and $\mathfrak h=\C$ in the above.  We obtain 
\begin{align*}
H^n(s\oplus\C,s\oplus\C)&=H^n(s\oplus\C,s)\oplus H^n(s\oplus\C,\C)\\
&=\bigoplus_{k+\ell=n}H^k(s,s)\tns H^\ell(\C,s)\oplus H^k(s,\C)\tns H^\ell(\C,\C)\\
&=\bigoplus_{k+\ell=n} H^k(s,\C)\tns H^\ell(\C,\C).
\end{align*}
Next, we need some elementary facts about the dimensions of $H^k(s,\C)$ and $H^\ell(\C,\C)$
\begin{align*}
&h^0(s,\C)=1,\quad h^1(s,\C)=0,\quad h^2(s,\C)=0,\quad h^3(s,\C)=1\\
&h^0(\C,\C)=1,\quad h^1(\C,\C)=1,\quad h^2(\C,\C)=0,\quad h^3(\C,C)=0.
\end{align*}
where $h^k=\dim(H^k)$.
Thus we obtain that 
\begin{align*}
h^0(s\oplus\C,s\oplus\C)=&h^0(s,\C)\cdot h^0(\C,\C)=1\cdot 1=1\\
h^1(s\oplus\C,s\oplus\C)=&h^1(s,\C)\cdot h^0(\C,\C)+h^0(s,\C)\cdot h^1(\C,\C)=0\cdot 1+1\cdot 1=1\\
h^2(s\oplus\C,s\oplus\C)=&h^2(s,\C)\cdot h^0(\C,C)+h^1(s,\C)\cdot h^1(\C,\C)+h^0(s,\C)\cdot h^2(\C,\C)\\
=&0\cdot 1+0\cdot 1+1\cdot 0=0\\
h^3(s\oplus\C,s\oplus\C)=&h^3(s,\C)\cdot h^0(\C,\C)+ h^2(s,\C)\cdot h^1(\C,\C)\\
&+h^1(s,\C)\cdot h^2(\C,\C)+h^0(s,\C)\cdot h^3(\C,\C)=1\cdot 1+0\cdot 1+0\cdot 0+1\cdot 1=1
\end{align*}
To compute $hc^n$, we get
\begin{align*}
&hc^0=z^1(s\oplus\C,\C)=1,\quad hc^1=h^2(s\oplus\C,\C)=0\\ &hc^2=h^3(s\oplus\C,\C)=1,\quad hc^3=h^4(s\oplus\C,\C)=1.
\end{align*}

Let us summarize these results in the following Table.

\begin{table}[ht]
  \begin{tabular}
    {lcccc}
    $n$&\vline&$HC^n$&$HRC^n$&$H^n$\\
    \hline
    0&\vline&1&1&1\\
    1&\vline&0&0&1\\
    2&\vline&1&0&0\\
    3&\vline&1&1&1\\
    \hline\\
  \end{tabular}
  \caption{Cohomology of the 4-dimensional Lie algebra $\sl(2,\C)\oplus\C$}
\end{table}

Just as in the case of the algebra $\sl(2,\C)$, we discover that there is a nontrivial cyclic deformation
$d_t=(1+t)d$, and we can use the identity transformation to see that this is formally equivalent to $d$, just as in the 3-dimesional case.

The real form $\sl(2,\R)\oplus\R$ has an invariant inner product which can be given by the same matrix as for the complex case, which has signature $(3,1)$. However, there is another real form with signature $(2,2)$.  For $\so(3,\R)\oplus\R$,
whose structure is given by $d=\psa{12}3-\psa{13}2+\psa{23}1$, we obtain an invariant inner product given by the matrix
$$B= \left[ \begin {array}{cccc} 1&0&0&0\\ \noalign{\medskip}0&1&0&0
\\ \noalign{\medskip}0&0&1&0\\ \noalign{\medskip}0&0&0&1\end {array}
 \right],
$$ whose signature is $(4,0)$. On the other hand, there is also a real form with signature $(3,1)$.

The reason that there are real forms with multiple signatures is that these algebras are direct sums of a simple and trivial algebra, and when combining the signatures, there are some variants possible.  For example, with $\sl(2,\R)$, the form of signature $(2,1)$ can combine with either the form of signature $(1,0)$ on $\R$ to give a form of signature $(3,1)$, or it can combine with the form of signature $(0,1)$ on $\R$ to give a form of signature $(2,2)$.  

The fact that there is an overlap in the possible signatures of these two real forms turns out to be important when we study deformations of the diamond and oscillator algebra.

The bases of the cohomology change from the $\sl(2,\R)\oplus\R$ case, but the dimensions of the cohomology are the same, so are also given by the
Table above. We have the same extra cyclic deformation for both real forms, realized by different cocycles. 
\subsection{The complex diamond algebra and its two real forms}

The complex diamond algebra can be given by $d=\psa{12}2-\psa{13}3+\psa{23}4$. It has an invariant inner product given by the matrix
$$B=
  \left[ \begin {array}{cccc} 0&0&0&1\\ \noalign{\medskip}0&0&1&0
\\ \noalign{\medskip}0&1&0&0\\ \noalign{\medskip}1&0&0&0\end {array}
 \right].
$$

The computation of the cohomology of the diamond algebra is not too difficult, but we omit it, for brevity, and merely summarize the results in the table below.

\begin{table}[ht]
  \begin{tabular}
    {lcccc}
    $n$&\vline&$HC^n$&$HRC^n$&$H^n$\\
    \hline
    0&\vline&1&1&1\\
    1&\vline&0&0&2\\
    2&\vline&1&1&2\\
    3&\vline&1&1&2\\
    \hline\\
  \end{tabular}
  \caption{Cohomology of the 4-dimensional complex diamond Lie algebra}

\end{table}
This is the first algebra we have encountered which has honest deformations, thus $H^2(d)$ does not vanish. In fact, a miniversal deformation of $d$ depends on two
parameters, say $t_1$ and $t_2$. We omit the formula for the
versal deformation, because it is long and we are interested in the cyclic versal deformation.

For the complex case, this just means it deforms to $\sl(2,\C)\oplus\C$. A basis for $HC^2$ is given by $\langle\psi=\psa{23}1-\psa{24}2+\psa{34}3\rangle$. 
The versal cyclic deformation is given by 
$$d_t=d+t\psi=\psa{12}2-\psa{13}3+\psa{23}4+t\psa{23}1-t\psa{24}2+t\psa{34}3.$$ It turns out that $d\sim d'=\sl(2,\C)\oplus \C$ with the same metric as for the diamond algebra.  In fact,
the matrix $G=\left[ \begin {array}{cccc} 0&0&1&1/2\\ \noalign{\medskip}1&0&0&0\\ \noalign{\medskip}0&{t}^{-1}&0&0\\ \noalign{\medskip}0&0&{t}^{-1}&-
1/2\,{t}^{-1}\end {array} \right] 
$
satisfies the property that $G(d_t(e_ie_j))=d'(G(e_i)G(e_j))$, which shows that $G$ gives an isomorphism betwee the algebras $d_t$ and $d'$.

\subsection{The real diamond algebra}
The structure we gave for the complex diamond Lie algebra coincides with the structure for the real diamond Lie algebra, so all of the above information on the
cohomology and cyclic cohomology is the same.  However, there is an important difference related to the fact that there are two real forms of the complex algebra.

In this case, we see that the signature of the real diamond algebra with the metric above is $(2,2)$, which is one of the possible signatures of $\sl(2,\R)\oplus\R$, but does not coincide with a possible signature for $\so(3)\oplus\R$, so it could not possibly deform to that algebra.  A simple computation shows that the real diamond algebra does deform to $\sl(2,\R)\oplus\R$.

The real diamond algebra is given as a semidirect product of the Heisenberg algebra by $\R$, and this fact plays a role in the applications of this algebra.

\subsection{The oscillator algebra} The oscillator algebra is the other real form of the complex diamond algebra.
The structure of this real Lie algebra can be given as 
$$d=\psa{12}3-\psa{13}2-\psa{23}4.$$ 
An invariant inner product is given by
$$\left[ \begin {array}{cccc} 0&0&0&1\\ \noalign{\medskip}0&1&0&0\\ \noalign{\medskip}0&0&1&0\\ \noalign{\medskip}1&0&0&0\end {array}
 \right].$$
 The signature of this form is $(3,1)$,  and the importance of this fact will become clear shortly.
 
 A basis for $HC^2$ is given by $\langle\psi=\psa{23}1+\psa{34}2-\psa{24}3\rangle$, so that once again, the versal deformation is of the form $d_t=d+t\psi$. It can be shown that this deformation is isomorphic to $\so(3,\R)\oplus\R$ when $t>0$ and to $\sl(2,\R)\oplus\R$ when $t<0$. It is important to note that both of these algebras have invariant inner products of signature $(3,1)$, and this explains why it is possible to have this deformation into two algebras, when we didn't have the same pattern for the diamond algebra. 
 
 A matrix of an isomorphism between the $d_t$ and $d'=\so(3,\R)\oplus\R$ is given by
 $G=\left[ \begin {array}{cccc} 0&0&1/2&1/2\\ 1/\sqrt{2t}&0&0\\ -1/\sqrt{2t}&0&0&0\\ \noalign{\medskip}0&0&1/(2t)&-1/
(2t)\end {array} \right] 
$. Note that $t$ must be positive for this matrix to be real.  We omit the transformation that gives the isomorphism between $d_t$ and $\sl(2,\R)\oplus\R$.

Now, a natural question arises about why we did not encounter the same issue about cyclic and reduced cyclic cohomology as for the cases of $\sl(2,\C)$ and $\sl(2,\C)\oplus\C$. Note that in those previous cases,  the algebra was not a coboundary in the cyclic cohomology, because the identity is not a cyclic cochain. This is also true in this case, but our algebra 
$d$ is still a coboundary because $[d,\ph^1_1-\ph^4_4]=d$.  The cochain $\ph^1_1-\ph^4_4$ is a cyclic cochain,  which means that $d$ is the coboundary of a cyclic cochain!

\section{5-dimensional Lie algebras with invariant inner product}

\subsection{The direct sum $\sl(2,\C)\oplus \C^2$ and its real forms}
This algebra can be given by $d=\psa{12}3+2\psa{23}2-2\psa{13}1$, and in this form, an invariant inner product is given by the matrix
$$B= \left[ \begin {array}{ccccc} 0&1&0&0&0\\ \noalign{\medskip}1&0&0&0&0\\ \noalign{\medskip}0&0&2&0&0\\ \noalign{\medskip}0&0&0&1&0
\\ \noalign{\medskip}0&0&0&0&1\end {array} \right].
$$
In the Table below, we summarize the cohomology information.  We omit any of the calculations used to obtain this information.

\begin{table}[ht]
  \begin{tabular}
    {lcccc}
    $n$&\vline&$HC^n$&$HRC^n$&$H^n$\\
    \hline
    0&\vline&2&2&2\\
    1&\vline&1&1&4\\
    2&\vline&1&0&2\\
    3&\vline&2&2&2\\
    \hline\\
  \end{tabular}
  \caption{Cohomology of the 5-dimensional Lie algebra $\sl(2,\C)\oplus\C^2$}
\end{table}
As usual for the algebras with a simple part, there is a metric deformation along the algebra itself, and the nontrivial cyclic 2-cocycle is just the cochain $d$ representing $\sl(2,\C)$.
As in the previous cases of $\sl(2,\C)$ and $\sl(2,\C)\oplus\C$, the cocycle $d$ is not a coboundary of a cyclic 1-cochain, but we still have $\exp(tI)d=\exp(t)d$, so we can eliminate the nontrivial cocycle by exponentiation of a 1-cochain that preserves the cyclic cocycles.

The real forms corresponding to this complex algebra are $\sl(2,\R)\oplus\R^2$ and $\so(3,\R)\oplus\R^2$. The cohomology dimensions remain the same, although the cocycles representing the cohomology and deformations change for the different algebras. 

As usual for the algebras with a simple part, there is a metric deformation along the algebra itself, and the nontrivial cyclic 2-cocycle is just $\sl(2,\C)$. 

The real forms corresponding to this complex algebra are $\sl(2,\R)\oplus\R^2$ and $\so(3,\R)\oplus\R^2$. The cohomology dimensions remain the same, although the cocycles representing the cohomology and deformations change for the different algebras. 

Since $\sl(2,\R)$ has an invariant metric of signature $(2,1)$, and since the signature of an invariant metric on $\R^2$ can be $(2,0)$ or $(1,1)$, this means we can obtain invariant metrics on $\sl(2,\R)\oplus\R^2$ of signature $(4,1)$, or $(3,2)$.

Similarly,  since $\so(3,\R)$ has an invariant metric of signature $(3,0)$, we can have invariant metrics of signature $(5,0)$ or $(4,1)$ on 
$\so(3,\R)\oplus\R^2$.  Note that there is an overlap in the possible signatures of the two real forms.

\subsection{The complex diamond algebra plus $\C$}.

The direct sum of the complex diamond algebra and $\C$ can be given by the structure $$d=\psa{12}2-\psa{13}3+\psa{23}4.$$
It has an invariant inner product given by the matrix
$$B=
\left[ \begin {array}{ccccc} 0&0&0&1&0\\ \noalign{\medskip}0&0&1&0&0\\ \noalign{\medskip}0&1&0&0&0\\ \noalign{\medskip}1&0&0&0&0
\\ \noalign{\medskip}0&0&0&0&1\end {array} \right].
$$

The versal cyclic deformation of the complex diamond algebra plus $\C$ can be given by
$$
d^\infty=\psa{12}2-\psa{13}3+\psa{23}4+t\psa{23}1-t\psa{24}2+t\psa{34}3.
$$
This deformation is isomorphic to $\sl(2,\C)\oplus\C^2$. Note that since the diamond algebra itself deforms to $\sl(2,\C)\oplus\C$,  this deformation is really obtained by just adding a plus $\C$ to each of the algebras.

We summarize the cohomology information in the Table below.

\begin{table}[ht]
  \begin{tabular}
    {lccccc}
    $n$&\vline&$HC^n$&$HRC^n$&$H^n$\\
    \hline
    0&\vline&2&2&2\\
    1&\vline&1&1&5\\
    2&\vline&1&1&5\\
    3&\vline&2&2&5\\
    \hline\\
  \end{tabular}
  \caption{Cohomology of the 4-dimensional complex diamond Lie algebra plus $\C$}
\end{table}

The deformation we gave for the complex diamond algebra plus $\C$ has real coefficients, so the structure also represents the cyclic versal deformation of $\sl(2,\R)\oplus\R$ with respect to the inner product given by the matrix above, which has signature $(3,2)$.  Both  $\sl(2,\R)\oplus\R^2$ and $\so(3,\R)\oplus\R^2$ have metrics of this signature, but the deformation is only isomorphic to  $\sl(2,\R)\oplus\R^2$ because the matrices of the transformations which give isomorphisms with $\so(3,\R)\oplus\R^2$ all have unavoidable complex coefficients.

Let us see what happens if we choose the metric given by the matrix 

$$B=
\left[ \begin {array}{ccccc} 0&0&0&1&0\\ \noalign{\medskip}0&0&1&0&0\\ \noalign{\medskip}0&1&0&0&0\\ \noalign{\medskip}1&0&0&0&0
\\ \noalign{\medskip}0&0&0&0&1\end {array} \right].
$$
This matrix has signature $(2,3)$ which is related to a matrix of signature $(3,2)$ in the way we described before, because multiplying the matrix by $-1$ reverses the signature, but doesn't affect the deformations.  It might seem that still, the cyclic versal deformation with this matrix might be different, but it turns out that the generator of $HC^2$ is the same for both matrices, so their versal deformations can coincide.

\subsection{The oscillator algebra plus $\R$}
We can use the same structure 
$$d=\psa{12}3-\psa{13}2-\psa{23}4$$ 
for the oscillator algebra plus $\R$ as we used for the oscillator algebra.  A matrix of an invariant inner product of signature $(4,1)$ is
$$
 \left[ \begin {array}{ccccc} 0&0&0&1&0\\ \noalign{\medskip}0&1&0&0&0
\\ \noalign{\medskip}0&0&1&0&0\\ \noalign{\medskip}1&0&0&0&0
\\ \noalign{\medskip}0&0&0&0&1\end {array} \right] 
$$
We can give the versal cyclic deformation of $d$ by the structure 
$$
d_t=\psa{12}3-\psa{13}2-\psa{23}4+t\psa{23}1+t\psa{34}2-t\psa{24}3.
$$
We compute that $d_t\sim\sl(2,\R)\oplus\R^2$ when $t<0$ and 
$d_t\sim\so(3,\R)\oplus\R^2$ when $t>0$.  This pattern corresponds to the pattern we observed for the deformations of the oscillator algebra.

As in the case of the oscillator algebra, if we use a metric of signature $(3,2)$, we obtain the same versal deformation as for the one of signature $(4,1)$, so there is no difference in the deformation pattern depending on the choice of metric.

Note that we don't have cochain representing the algebra appearing as a nontrivial cocycle, because if we choose $\beta=\pha11-\pha44$, then $[\beta,d]=d$ so the algebra is given by a trivial cocycle.
\subsection{The algebra $W_3$}

The algebra $W_3$ can be given by the cochain $$d=\psa{34}2+\psa{35}1+\psa{45}3.$$ 
This is the first nilpotent metric Lie algebra we have encountered in this paper.

It has an invariant inner product given by the matrix
$$B=
\left[ \begin {array}{ccccc} 0&0&0&-1&0\\ \noalign{\medskip}0&0&0&0&1\\ \noalign{\medskip}0&0&1&0&0\\ \noalign{\medskip}-1&0&0&0&0
\\ \noalign{\medskip}0&1&0&0&1\end {array} \right] 
$$

The cohomology of the algebra $W_3$ is summarized in the Table below.

\begin{table}[ht]
  \begin{tabular}
    {lccccc}
    $n$&\vline&$HC^n$&$HRC^n$&$H^n$\\
    \hline
    0&\vline&2&2&2\\
    1&\vline&3&3&7\\
    2&\vline&3&3&9\\
    3&\vline&2&2&9\\
    \hline\\
  \end{tabular}
  \caption{Cohomology of the 5-dimensional complex algebra $W_3$}

\end{table}
A basis of the cyclic nontrivial 2-cocycles (that is, a basis of $HC^2$) is given by \begin{align*}HC^2=\langle&
\delta^1=-\psa{23}1+\psa{35}1-\psa{24}3+\psa{45}3+\psa{34}5\\
&\delta^2=\psa{23}2-\psa{35}2-\psa{25}3+\psa{35}5,\\
&\delta^3=\psa{13}2-\psa{15}3-\psa{35}4
\rangle\end{align*}

This is  the first algebra for which $HC^2$ is larger than 1-dimensional, and moreover, the versal deformation determined by the above basis of $H^2$ has higher order terms. In fact, a versal deformation 
$d^\infty$ can be given by 
$$
d^\infty=d+t_1\delta^1+t_2\delta^2+t_3\delta^3
+t_1t_3(-\psa{15}3-\psa{35}4-\psa{12}3+\psa{23}4+\psa{13}5).
$$
There are no relations on the base, so that $[d^\infty,d^\infty]=0$,
and thus $\dinfty$ is a Lie algebra for all values of the parameters 
$t_i$.

On the surface
$$
2t_1t_3(2+t_1)=t_2^2
$$
we have that $\dinfty$ is isomorphic to the diamond algebra plus $\C$,
while elsewhere it is isomorphic to $\sl(2,\C)\oplus \C^2$. This is a typical pattern in the deformation picture of algebras.

For example, $d_t=d+t_2\delta^2$ is isomorphic to $\sl(2,\C)\oplus\C^2$ when $t_2\ne0$, while
the 1-parameter deformations $d_t=d+t_1\delta^1$ and $d_t=d+t_3\delta^3$ both give jump deformations to the diamond algebra plus $\C$. Notice that we obtain 2 different 1-parameter families isomorphic to the same algebra, but the generic deformation is still to the algebra $\sl(2,\C)\oplus\C^2$.

\section{6-dimensional Lie algebras with invariant inner product}

Up to now, our method of determining the Lie algebras with invariant inner product depended on the classifications of the moduli spaces of algebras of the corresponding dimension given in 
\cite{fp3,Ott-Pen,fp8,fp20}, which built upon earlier classifications of the algebras.  For dimension 6, we have not constructed the moduli space of such algebras, so there are two methods of attack which will yield all the metric Lie algebras.

First, there is a classification due to Mubarakzyanov \cite{mub3}, who gave a classification of all solvable, nonnilpotent real Lie algebras of dimension 6.  Since the classification of nilpotent algebras is known \cite{FaSa}, and the classification of nonsolvable Lie algebras is straightforward, this is sufficient to give a complete classification of the metric Lie algebras.

The second method is due to Kac \cite{kac}, who gave a method of constructing a metric Lie algebra of dimension $n+2$, given any metric Lie algebra of dimension $n$. This construction yields an algebra with a nontrivial center, but from \cite{FaSa}, we know that any solvable metric Lie algebra has a nontrivial center. Moreover, this construction yields every solvable metric Lie algebra, as well as some of the nonsolvable Lie algebras.

We first study the algebras which are nonsolvable, then the algebras like the diamond algebra plus $\C^2$, $W_3\oplus\C$ and the algebra $W_4$,  the latter being the only nilpotent 6-dimensional metric Lie algebra, and then we will study the rest of the solvable metric Lie algebras.

 \subsection{The direct sum $\sl(2,\C)\oplus \sl(2,\C)$}
 The algebra can be given by 
 $$d=\psa{12}3+2\psa{23}2-2\psa{13}1+\psa{45}6+2\psa{56}5-2\psa{46}4.$$  It has an invariant inner product given by the matrix 
 $$B=\left[ \begin {array}{cccccc} 0&1&0&0&0&0\\ \noalign{\medskip}1&0&0&0&0&0\\ \noalign{\medskip}0&0&2&0&0&0\\ \noalign{\medskip}0&0&0&0&1&0
\\ \noalign{\medskip}0&0&0&1&0&0\\ \noalign{\medskip}0&0&0&0&0&2
\end {array} \right].
$$
We have $$HC^2=\langle \psi^1=-\psa{46}4+\psa{56}5+\tfrac12\psa{45}6,\quad
\psi^2=-\psa{13}1+\psa{23}2+\tfrac12\psa{12}3\rangle,$$
which are essentially the 2-cochains which determine the direct summands that make up the algebra.
The metric versal deformation is given by 
$$d^\infty=d+t_1\psi^1+t_2\psi^2,$$
and this is just the same type of pattern of deforming along itself that we saw with simple Lie algebras, except that now there two parameters of deformation along the original algebra, induced by the two simple pieces of the algebra. 
We summarize the cohomology of this algebra in the following table.
\begin{table}[ht]
  \begin{tabular}
    {lccccc}
    $n$&\vline&$HC^n$&$HRC^n$&$H^n$\\
    \hline
    0&\vline&0&0&0\\
    1&\vline&0&0&0\\
    2&\vline&2&0&0\\
    3&\vline&0&0&0\\
    \hline\\
  \end{tabular}
  \caption{Cohomology of the 6-dimensional complex algebra $\sl(2,\C)\oplus\sl(2,\C)$}
 \end{table} 
 
 Note that the only variation in the dimensions of $HC^n$ is in the case $n=2$, where we pick up the special deformation along the algebra itself.  It is not difficult to show that no nontrivial linear combination of the cocycles $\psi^1$ and $\psi^2$ is a coboundary of a cyclic 1-cochain.
 
 \subsection{The semidirect product $T^*(\sl(2,\C))$ of $\sl(2,\C)$ and $\C^3$}
  This is the algebra which is given by the coadjoint representation of $\sl(2,\C)$ on its cotangent space.
  This algebra is given by the cochain
  $$d=2\psa{12}2-2\psa{13}3+\psa{23}1+2\psa{14}4-2\psa{16}6+2\psa{25}4+\psa{26}5+\psa{34}5+2\psa{35}6.
  $$
  It has an invariant inner product given by the matrix
  $$B=\left[ \begin {array}{cccccc} 0&0&0&0&2&0\\ \noalign{\medskip}0&0&0&0&0&1\\ \noalign{\medskip}0&0&0&-1&0&0\\ \noalign{\medskip}0&0&-1&0&0&0
\\ \noalign{\medskip}2&0&0&0&0&0\\ \noalign{\medskip}0&1&0&0&0&0
\end {array} \right]. 
 $$
  We have 
  $$HC^2=\langle \psi^1=\tfrac12\psa{46}1-\psa{45}2+\psa{56}3,\psi^2=-2\psa{12}4+\psa{23}5-2\psa{13}6\rangle.
  $$
  
  The versal metric deformation is given by 
  \begin{align*}
  d^\infty=&d+t_1\psi^1+t_2\psi^2+t_1t_2(\cdots),
  \end{align*}
  where the $\cdots$ stands for a sum of other cochain terms.  There are some terms with a factor of $1+\sqrt{1+4t_1t_2}$ in the denominator, but these denominators are constant if either $t_1$ or $t_2$ vanish. 
  
  When $t_1\ne 0$ the deformation is isomorphic to $\sl(2,\C)\oplus\sl(2,\C)$. In particular
  the 1-parameter family $d_t=d+t_1\psi^1$ gives a jump deformation to this algebra. On the other hand, when $t_1=0$, the deformation is $d_t=d + t_2\psi^2$  is isomorphic to itself.  The explanation is different than the first case, becauses the cochain $d$ actually is a coboundary. This type of behavior occurs for ordinary deformations, but is a rare phenomena. 
  
  Since the algebra does not occur as a deformaton of itself in the ordinary deformation sense, we must have $\psi^2$ is a cyclic coboundary of a noncyclic 1-cochain.  In fact, the cochain $\beta=\pha{1}5-2\pha{2}4$ satisfies $[d,\beta]=\psi^2$, and it is not difficult to see that $\beta$ is not a cyclic cochain.
  
   This example makes it clear that an algebra which has a simple quotient can have a higher dimensional $HC^2$ even when it is not given by a direct summand.  
  
  The following Table gives the cohomology of this algebra in low degrees.
  \begin{table}[ht]
  \begin{tabular}
    {lccccc}
    $n$&\vline&$HC^n$&$HRC^n$&$H^n$\\
    \hline
    0&\vline&0&0&0\\
    1&\vline&0&0&1\\
    2&\vline&2&1&1\\
    3&\vline&0&0&0\\
    \hline\\
  \end{tabular}
  \caption{Cohomology of the 6-dimensional complex algebra $T^*(\sl(2,\C))$}
 \end{table} 
  
 \subsection{The direct sum $\sl(2,\C)\oplus\C^3$}
 For this algebra, defined by $$d=\psa{12}3-2\psa{13}1+2\psa{23}2$$ has an invariant inner product given by the matrix 
$$B=\left[ \begin {array}{cccccc} 0&1&0&0&0&0\\ \noalign{\medskip}1&0&0&0&0&0\\ \noalign{\medskip}0&0&2&0&0&0\\ \noalign{\medskip}0&0&0&1&0&0
\\ \noalign{\medskip}0&0&0&0&1&0\\ \noalign{\medskip}0&0&0&0&0&1
\end {array} \right] 
$$
 $$HC^2=\langle \psi^1=\psa{56}4-\psa{46}5+\psa{45}6,\quad\psi^2=-\psa{13}1+\psa{23}2+\tfrac12\psa{12}3\rangle.
  $$
  The metric versal deformation of this algebra is just
  $$d^\infty=d+t_1\psi^1+t_2\psi^2$$
  and there are no relations on the base.  The algebra deforms to $\sl(2,\C)\oplus\sl(2,\C)$ as well as to itself, in the usual pattern for algebras with a simple quotient. Note that the second cocycle is just half of $d$. When $t_1\ne0$, the deformation is isomorphic to $\sl(2\C)\oplus\sl(2,\C)$, while if $t_1=0$, the deformation is isomorphic to $d$.

We summarize the cohomology of this algebra in the following Table.
\begin{table}[ht]
  \begin{tabular}{lccccc}
    $n$&\vline&$HC^n$&$HRC^n$&$H^n$\\
    \hline
    0&\vline&3&3&3\\
    1&\vline&3&3&9\\
    2&\vline&2&1&9\\
    3&\vline&3&3&6\\
    \hline\\
  \end{tabular}
  \caption{Cohomology of the 6-dimensional complex algebra $\sl(2,\C)\oplus\C^3$}
 \end{table}

 \subsection{The algebra family \oscl}
 The algebra \oscl (see \cite{pdl}) is represented by 
 $$
 d=\psa{12}2+\lambda\psa{13}3-\psa{14}4-\lambda\psa{15}5+\lambda\psa{35}6+\psa{24}6,
 $$
 where the case $\lambda=0$ is generally excluded from the family, because it is decomposable. 

 This algebra has some symmetries in the sense that replacing $\lambda$ with $1/\lambda$ (when $\lambda\ne 0$) or $-\lambda$ gives an isomorphic algebra.
 
 A matrix of an invariant bilinear form for this algebra is
 $$
 B=\left[ \begin {array}{cccccc} 0&0&0&0&0&1\\ \noalign{\medskip}0&0&0&1&0&0\\ \noalign{\medskip}0&0&0&0&1&0\\ \noalign{\medskip}0&1&0&0&0&0
\\ \noalign{\medskip}0&0&1&0&0&0\\ \noalign{\medskip}1&0&0&0&0&0
\end {array} \right] .
$$
Note that since $\lambda$ does not appear in the matrix, the same bilinear form works for all $\lambda$.

In addition to $\lambda=0$, the case $\lambda=1$ is special in the sense that for both cases, the cohomology and deformation patterns are not generic. 

We first give the generic picture of the deformations. We have
$$HC^2=\langle\psi^1=-\psa{24}1+\psa{26}2-\psa{46}4+\lambda\psa{35}1-\lambda\psa{36}3+\lambda\psa{56}5,\quad \psi^2=\psa{12}2-\psa{14}4+\psa{24}6\rangle.
$$
This basis works unless $\lambda=0$, which we will study separately, and $\lambda=1$, which has much larger cohomology. The reason the basis is not valid for $\lambda=0$ is that the first cocycle was obtained from a cocycle with $\lambda$ in the denominator, which is invalid if $\lambda=0$.

The versal deformation is given by 
$$d^\infty=d+t_1\psi^1+t_2\psi^2+t_1t_2(-\psa{24}1+\psa{26}2-\psa{46}4).
$$
When $t_1\ne0$, the deformation is equivalent to $\sl(2,\C)\oplus\sl(2,\C)$.
When $t_1=0$, the deformation is a smooth deformation in a neighborhood of $\mathfrak{g}_{6,2}(\lambda)$, and these are all the deformations that occur generically. The deformation is isomorphic to $\mathfrak{g}_{6,2}(\tfrac\lambda{1+t_2})$. This isomorphism can be given by the linear transformation given by the  matrix 
$$ \left[ \begin {array}{cccccc} 1&0&0&0&0&0\\\noalign{\medskip}0&0&1&0&0
&0\\\noalign{\medskip}0&1&0&0&0&0\\\noalign{\medskip}0&0&0&0&1&0
\\\noalign{\medskip}0&0&0&1&0&0\\\noalign{\medskip}0&0&0&0&0&1
\end {array} \right].
$$
 
 The generic cohomology picture is given in the Table below.
 \begin{table}[ht]
  \begin{tabular}
    {lccccc}
    $n$&\vline&$HC^n$&$HRC^n$&$H^n$\\
    \hline
    0&\vline&1&1&1\\
    1&\vline&1&1&5\\
    2&\vline&2&2&7\\
    3&\vline&1&1&6\\
    \hline\\
  \end{tabular}
  \caption{Cohomology of the algebra \oscl  for generic $\lambda$}
  \end{table}
 
\subsubsection{The special case $\lambda=0$}

This is isomorphic to the diamond algebra plus $\C^2$, which can be represented by the 
cochain 
$$d=\psa{12}3-\psa{13}2+\psa{23}4.$$
The algebra $\mathfrak{g}_{6,2}(0)$ is represented by the cochain 
$$\psa{12}2-\psa{14}4+\psa{24}6.
$$
It has the same bilinear form as the generic case, so we will use the $\mathfrak{g}_{6,2}(0)$  structure in our computations below.

We do have to replace the generic cocycles representing $H^2$, and we have
$$HC^2=\langle \psi^1=\psa{24}1-\psa{26}2+\psa{46}4,\quad\psi^2=\psa{13}3-\psa{15}5+\psa{35}6\rangle.
$$
The versal deformation is given by 
$$d^\infty=d+t_1\psi^1+t_2\psi^2+t_1t_2(-\psa{35}1+\psa{36}3-\psa{56}5).$$

When $t_1\ne0$ and $t_2\ne0$, it jumps to $\sl(2,C)\oplus\sl(2,\C)$. For example,
the deformation given by 
$$d_t=d+t(\psi^1+\psi^2)+t^2(-\psa{35}1+\psa{36}3-\psa{56}5),
$$
which is obtained by substituting $t_1=t$ and $t_2=t$, is isomorphic to $\sl(2,\C)\oplus\sl(2,\C)$ whenever $t\ne0$. A matrix of a linear transformation which realizes this isomorphism is 
$$
\left[ \begin {array}{cccccc} 0&0&0&t_{{1}}&0&0\\\noalign{\medskip}0&1&0&0&0&0\\\noalign{\medskip}-1/2&0&0&0&0&-1/2\,t_{{1}}
\\\noalign{\medskip}0&0&-{t_{{1}}}^{3}&0&0&0\\\noalign{\medskip}0&0&0&0
&1&0\\\noalign{\medskip}1/2\,t_{{1}}&0&0&0&0&-1/2\,{t_{{1}}}^{2}
\end {array} \right].
$$
The deformation $d_t=d+t_1\psi^1$ is isomorphic to $\sl(2,\C)\oplus\C^3$. Finally, if we substitute $t_1=0$ into $d^\infty$, we obtain the algebra $\mathfrak{g_{6,2}}(t_2)$.  This means that when $t_2$ is small, we obtain an element of the family which is ``near" to our algebra. This is the simplest type of an example of a smooth deformation. 

\begin{table}[ht]
  \begin{tabular}
    {lccccc}
    $n$&\vline&$HC^n$&$HRC^n$&$H^n$\\
    \hline
    0&\vline&3&3&3\\
    1&\vline&3&3&10\\
    2&\vline&2&2&13\\
    3&\vline&3&3&12\\
    \hline\\
  \end{tabular}
  \caption{Cohomology of the complex diamond algebra plus $\C^2$}
 \end{table} 
  \subsubsection{The special case $\lambda=1$}
  
  The first evidence that this case is special is that the dimension of $HC^2$ is 6, rather than the generic value 2.
  We have 
 \begin{align*}
 HC^2=\langle&
 \psi^1=\psa{13}3-\psa{15}5+\psa{35}6,\quad
 \psi^2=\psa{34}1-\psa{36}2+\psa{46}5\\&
 \psi^3=\psa{25}1-\psa{26}3+\psa{56}4,
 \quad
 \psi^4=-\psa{24}1+\psa{26}2-\psa{46}4+\psa{35}1-\psa{36}3+\psa{56}5\\&
 \psi^5=\psa{12}3-\psa{15}4+\psa{25}6,\quad
 \psi^6=\psa{13}2-\psa{14}5+\psa{34}6\rangle.
 \end{align*}
 We get
 $$d^\infty=d+t_i\psi^i+
 \frac{t_1t_4+t_2t_5+t_3t_6}{1+t_1}(\psa{24}1+\psa{26}2+\psa{46}4).
 $$
 Also, there are 3 relations on the base, given by 
 \begin{align*}
     -2\,t_{{2}}t_{{4}}-4\,t_{{5}}t_{{3}}-2\,t_{{5}}t_{{4}}t_{{3}}+2\,{t_{{
5}}}^{2}t_{{1}}+2\,t_{{5}}t_{{6}}t_{{2}}&\\
2\,t_{{1}}t_{{4}}+4\,t_{{6}}t_{{3}}+2\,t_{{6}}t_{{4}}t_{{3}}-2\,t_{{6}
}t_{{5}}t_{{1}}-2\,{t_{{6}}}^{2}t_{{2}}&\\
2\,t_{{1}}t_{{5}}-2\,t_{{6}}t_{{2}}&
 \end{align*}
 These relations are polynomials representing the coefficients of the three basis elements of $H^3$ in the bracket $[\dinfty,\dinfty]$ (these are the only terms in the bracket), and therefore these relations must be equal to zero.
 When solving these relations in Maple, we came up with 6 solutions.  Each of the solutions is local, meaning that the origin in $t$ space lies within the range of the solution.  We do not give these solutions, but merely give the results of the analysis of the deformations given by the solutions.  
 
 The algebra $\mathfrak g_{6,2}(1)$ has jump deformations to the algebras $\sl(2,\C)\oplus\sl(2,\C)$
 $\sl(2,\C) \oplus \C^3$,
 $T^*(\sl(2,\C))$ and $\mathfrak g_{6,3}$, the algebra we will study next, as well as smooth deformations in a neighborhood  of itself in \oscl.  The value $\lambda=1$ is truly a special case. Later we will discuss an interpretation which will make this case even more special. 
 
  The cohomology of $\mathfrak g_{6,2}(1)$ is given in the table below.
   \begin{table}[ht]
  \begin{tabular}
    {lccccc}
    $n$&\vline&$HC^n$&$HRC^n$&$H^n$\\
    \hline
    0&\vline&1&1&1\\
    1&\vline&3&3&5\\
    2&\vline&6&6&7\\
    3&\vline&3&3&6\\
    \hline\\
  \end{tabular}
  \caption{Cohomology of the algebra $\mathfrak g_{6,2}(1)$}
  \end{table}
  
  \subsection{The algebra $\mathfrak{g}_{6,3}$} The algebra $\mathfrak{g}_{6,3}$ is represented by the cochain
  $$
  d=\psa{34}4+\psa{35}4+\psa{35}5+\psa{13}1+\psa{13}2+\psa{23}2-\psa{14}6-\psa{15}6-\psa{25}6.
  $$
  A matrix of an invariant bilinear form for the algebra is
  $$
  B=\left[ \begin {array}{cccccc} 0&0&0&1&0&0\\ \noalign{\medskip}0&0&0&0&1&0\\ \noalign{\medskip}0&0&0&0&0&1\\ \noalign{\medskip}1&0&0&0&0&0
\\ \noalign{\medskip}0&1&0&0&0&0\\ \noalign{\medskip}0&0&1&0&0&0
\end {array} \right] 
.$$
A basis for $HC^2$ is given by 
$$
HC^2=\langle \psi^1=w-\psa{16}2+\psa{15}3+\psa{56}4,\quad\psi^2=\psa{23}1+\psa{34}5-\psa{24}6\rangle.
$$

  The versal metric deformation is given by 
  \begin{align*}
  d^\infty=&d+t_1\psi^1+t_2\psi^2+\ho,
  \end{align*}
  where all the higher order terms are multiplied by $t_1t_2$, so they are only quadratic in the parameters.  
  
  The algebra has jump deformations  to $\sl(2,\C)\oplus\sl(2,\C)$ and to 
  $T^*(\sl(2,\C))$. It also deforms in a neighborhood of $\mathfrak{g}_{6,2}(1)$ but does not have a jump deformation to that algebra. For example, if both $t_1$ and $t_2$ are nonzero, then the deformation is isomorphic to $\sl(2,\C)\oplus\sl(2,\C)$, while if $t_2=0$, then it gives a jump deformation to $T^*(\sl(2,\C))$.   Finally, if $t_1=0$, then we get a deformation in a neighborhood of $\mathfrak{g}_{6,2}(1)$, but it does not jump to $\mathfrak{g}_{6,2}(1)$.
  
  This means that there are two distinct algebras which deform in a neighborhood of $\mathfrak{g}_{6,2}(1)$. In some prior works (see \cite{fp16}), the authors have encountered this type of situation, and the analysis sometimes required an interchange of the elements in the family to fit the deformation pattern. After the completion of the description of the deformations of the algebras, we will address this situation.

  We summarize the cohomology of this algebra in the following Table.
\begin{table}[ht]
  \begin{tabular}
    {lccccc}
    $n$&\vline&$HC^n$&$HRC^n$&$H^n$\\
    \hline
    0&\vline&1&1&1\\
    1&\vline&1&1&3\\
    2&\vline&2&2&3\\
    3&\vline&1&1&2\\
    \hline\\
  \end{tabular}
  \caption{Cohomology of the algebra $\mathfrak{g}_{6,3}$}
 \end{table}

\subsection{The algebra $W_3\oplus\C$}

The algebra is given by $$d=\psa{34}2+\psa{35}1+\psa{45}3.$$
  It has an invariant matrix given by
  $$B=\left[ \begin {array}{cccccc} 0&0&0&-1&0&0\\ \noalign{\medskip}0&0&0&0&1&0\\ \noalign{\medskip}0&0&1&0&0&0\\ \noalign{\medskip}-1&0&0&0&0&0
\\ \noalign{\medskip}0&1&0&0&0&0\\ \noalign{\medskip}0&0&0&0&0&1
\end {array} \right].$$

A basis for the cohomology $HC^2$ is given by 
\begin{align*}HC^2=\langle& \psi^1=\psa{26}1+\psa{46}5+\psa{24}6,\psi^2=\psa{16}2+\psa{56}4-\psa{15}6\\&
\psi^3=\psa{16}1-\psa{46}4+\psa{14}6-\psa{26}2+\psa{56}5+\psa{25}6,
\psi^4=-\psa{23}1-\psa{24}3+\psa{34}5,\\&
\psi^5=-\psa{13}2+\psa{15}3+\psa{35}4\\
\psi^6=\psa{23}2-\psa{25}3+\psa{35}5
\rangle
\end{align*}
We have 
$$d^\infty=d+t_i\psi^i+\ho,$$ where $\ho$ stands for the higher order terms all of which have degree 2 in the parameters.

There are 3 relations on the base, 2 of which are nontrivial. 
We have jump deformations to $\sl(2,\C)\oplus\sl(2,\C)$, $T^*(\sl(2,\C))$, $\sl(2,\C)\oplus\C^3$, and $\mathfrak{g}_{6,3}$, and smooth deformations in a neighborhood of $\mathfrak{g}_{6,2}(1)$.  Note that it does not have a jump deformation to $\mathfrak{g}_{6,2}(1)$, which is very important.

We summarize the cohomology of this algebra in the Table below.
\begin{table}[ht]
  \begin{tabular}
    {lccccc}
    $n$&\vline&$HC^n$&$HRC^n$&$H^n$\\
    \hline
    0&\vline&3&3&3\\
    1&\vline&5&5&12\\
    2&\vline&6&6&21\\
    3&\vline&5&5&24\\
    \hline\\
  \end{tabular}
  \caption{Cohomology of the complex algebra $W_3\oplus \C$.}
 \end{table}

 \subsection{The complex algebra $W_4$} 
 
 The nilpotent algebra $W_4$ is given by $$d=\psa{15}4-\psa{13}2-\psa{35}6.$$  It has an invariant inner product given by 
 $$B=\left[ \begin {array}{cccccc} 0&0&0&0&0&1\\ \noalign{\medskip}0&0&0&0&1&0\\ \noalign{\medskip}0&0&0&1&0&0\\ \noalign{\medskip}0&0&1&0&0&0
\\ \noalign{\medskip}0&1&0&0&0&0\\ \noalign{\medskip}1&0&0&0&0&0
\end {array} \right] 
 $$
 We have 
 \begin{align*}
     HC^2=&\langle
 \psi^1=\psa{34}1-\psa{36}3+\psa{46}4+\psa{25}1-+\psa{26}2+\psa{56}6,\quad
 \psi^2=\psa{15}1-\psa{26}2+\psa{56}6\\&
 \psi^3=\psa{24}2-\psa{25}3+\psa{45}5+\psa{14}4-\psa{16}3+\psa{46}6,\quad
 \psi^4=\psa{45}1-\psa{46}2+\psa{56}3,\\&
 \psi^5=\psa{35}1-\psa{36}2+\psa{56}4,\quad
 \psi^6=\psa{23z}2-\psa{25}4+\psa{35}5,\quad
 \psi^7=\psa{23}1-\psa{26}4+\psa{36}5,\\&
 \psi^8=\psa{14}2-\psa{15}3+\psa{45}6,\quad
 \psi^9=\psa{12}4-\psa{13}5+\psa{45}6+\psa{23}6,\quad
 \psi^{10}=\psa{12}3-\psa{14}5+\psa{24}6\\&
 \psi^{11}=\psa{12}2-\psa{15}5+\psa{25}6,\quad
 \psi^{12}=\psa{12}1-\psa{16}5+\psa{26}2+\psa{23}3-\psa{24}4+\psa{34}5
 \rangle
 \end{align*}
 The form of the versal deformation is relatively simple:
 $$d^\infty=d+t_i\psi^i+\ho,
 $$
 where the higher order terms are all quadratic or cubic in the parameters.
 
 The complexity of this versal deformation is in terms of the relations on the base.  The 8 relations give rise to 21 different solutions to the relations, yielding jump deformations to every 6-dimensional metric Lie algebra other than itself,  including \oscl    for all values of 
 $\lambda$.
 
 We summarize the cohomology information for $W_4$ in the Table below.
\begin{table}[ht]
  \begin{tabular}
    {lccccc}
    $n$&\vline&$HC^n$&$HRC^n$&$H^n$\\
    \hline
    0&\vline&3&3&3\\
    1&\vline&8&8&15\\
    2&\vline&12&12&30\\
    3&\vline&8&8&36\\
    \hline\\
  \end{tabular}
  \caption{Cohomology of the complex algebra $W_4$.}
 \end{table}
 
 \section{Ordering of the elements in dimension 6}
 
It is natural to take $\sl(2,\C)\oplus\sl(2,\C)$ to be the first element, or algebra 1, while $T^*(\sl(2,\C))$ is 2, and $\sl(2,\C)\oplus\C^3$ is algebra 3. 

There are some problems involving the algebras $\mathfrak g_{6,2}(1)$ and $\mathfrak g_{6,3}$. They arise because of the deformation patterns involving $g_{6,2}(1)$. Both this algebra and $\mathfrak g_{6,3}$ deform in a neighborhood of $g_{6,2}(1)$. Also, the cohomology of $\mathfrak g_{6,2}(1)$ has dimension $6$, while that of $\mathfrak g_{6,3}$ has only dimension 2. In \cite{fp16}, we found that in this circumstance, it makes sense to interchange the two elements,
in other words,  we have a family $4(\lambda)$, in which for $\lambda\ne1$, the algebra is $\mathfrak g_{6,2}(\lambda)$, but $4(1)=\mathfrak g_{6,3}$. Then the algebra $\mathfrak g_{6,2}(1)$ becomes the algebra 5.  It fits  the deformation picture in that we have algebras with a higher cohomology dimension jump to algebras with a lower cohomology dimension.  But there is another way in which this switch fits the deformation pattern.  The algebra $W_3\oplus\C$ deforms in a neighborhood of $\mathfrak g_{6,2}(1)$ but does not jump to that algebra.  With the switch, we obtain that 
$W_3 \oplus\C$, to which we assign the algebra number 6, jumps to $4(1)$ and deforms in a neighborhood of that algebra.  

Finally, we assign the number 7 to the algebra $W_4$.

\vskip .1in
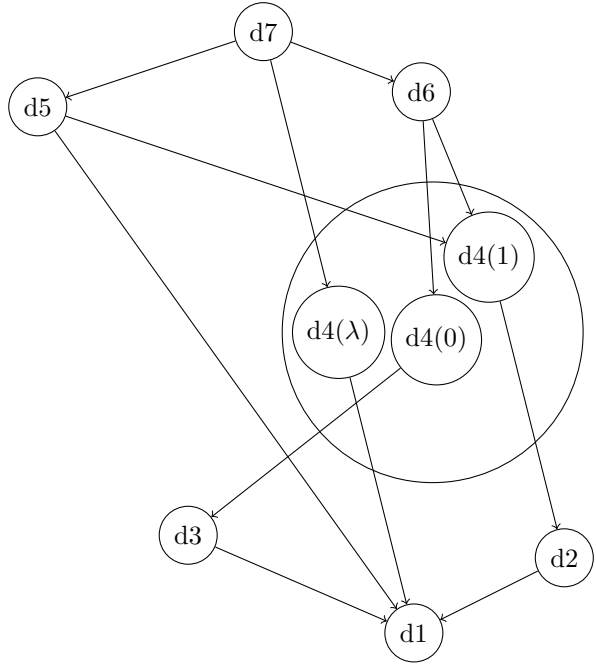
\begin{center}
\begin{tikzpicture}
\node [circle, draw] at (0,6) (d7) {d7};
\node [circle, draw] at (2.1,5.2) (d6) {d6};
\node [circle, draw] at (-3,5) (d5) {d5};

\node[] at (2.25,2) (center) {};
\draw (center) circle (2);
\node [circle, draw] at (3,3) (d41) {d4(1)};
\node [circle, draw] at (2.3,1.9) (d40) {d4(0)};
\node [circle, draw] at (1, 2) (d4lambda) {d4($\lambda$)};

\node [circle, draw] at (-1,-.7) (d3) {d3};
\node [circle, draw] at (4,-1) (d2) {d2};
\node [circle, draw] at (2,-2) (d1) {d1};

\draw[->] (d7)--(d6);

\draw[->] (d7)--(d5);
\draw[->] (d6)--(d41);
\draw[->] (d5)--(d41);
\draw[->] (d41)--(d2);
\draw[->] (d2)--(d1);
\draw[->] (d4lambda)--(d1);
\draw[->] (d40)--(d3);
\draw[->] (d3)--(d1);
\draw[->](d6)--(d40);
\draw[->](d5)--(d1);
\draw[->](d7)--(d4lambda);

\end{tikzpicture}
\end{center}
\vskip .01in

In the picture, the downward direction corresponds to the direction of the jump deformation.  The fact that it is always possible to arrange the ordering of the algebras to have this property is from a standard fact of deformation theory: every jump deformation leads to an algebra with strictly smaller second cohomology group, which means there is a unidirectionality in the jump deformation theory. For smooth deformations along a family, the dimension of the cohomology need not drop, and there are often special points in the family where the cohomology has larger dimension than the generic case.  In our picture, the corresponding points are $d4(1)$ and $d4(0)$.
\newpage
The Table below gives the same information as the picture with some more details.
\begin{table}[ht]
\begin{tabular}
{llcc}
Type&Name&Jump Deformations&Smooth Deformations\\
\hline
1&$\sl(2,\C)\oplus\sl(2,\C)$&&\\
2&$T^*(\sl(2,\C))$&1&\\
3&$\sl(2,\C)\oplus\C^3$&1&\\
$4(\lambda), \lambda \neq 0,1$&$\mathfrak g_{6,2}(\lambda)$&1&$4(\lambda)$\\
$4(0)$&$\mathfrak g_{6,2}(0)=\text{diamond}\oplus\C^2$&1,3&$4(0)$\\
$4(1)$&$\mathfrak g_{6,3}$&1,2&$5$\\
5&$\mathfrak g_{6,2}(1)$&1,2,3,$4(1)$&$5$\\
6&$W_3\oplus\C$&1,2,3,4(0),4(1)&5 \\
7&$W_4$&1,2,3,4,5,6&\\  
\end{tabular}
\end{table}

\newpage
 \section{Conclusions}
 
In our investigations of metric algebras of dimension up to 6, exactly the algebras with a simple quotient are the ones for which the reduced cyclic cohomology has smaller dimension than the cyclic cohomology, and these algebras have smooth deformations along themselves.  Whether this pattern holds in general is not obvious, because we have only looked at low dimensional algebras, and  the only example of an algebra with a simple quotient but not a simple direct summand that we encountered is 
$T^*(\sl(2,\C)$, which is a semidirect product of $\sl(2,\C)$ and $\C^3$. 

In all of our examples, we have seen that the reduced cyclic cohomology coincides with the reduced cyclic cohomology except in the cases where there is a simple quotient, and those algebras have a formal nontrivial deformation along themselves.  It is reasonable to conjecture that this is the general pattern.

For the algebra $\sl(2,\C)$, we found that the cochain $d$ representing the algebra is not a coboundary, which never happens for ordinary deformations, as it is the coboundary of the identity transformation.  For cyclic deformations, the problem begins with the fact that the identity is never a cyclic 1-cochain for a simple algebra, because of the nature of the Killing form. However, we defined a notion of reduced cyclic cohomology, which considers as coboundaries those cochains which are cyclic coboundaries, rather than coboundaries of cyclic cochains. This increases the dimension of the cyclic coboundaries, and so decreases the dimension of the cyclic cohomology, which we called reduced cyclic cohomology.

In addition, we noted that the bracket of the identity matrix with any cochain is itself, 
therefore applying the exponential of a multiple of the identity to a cyclic deformation gives another cyclic deformation. If we consider any two such deformations to be equivalent, then
in this sense, replacing the algebra with a nonzero multiple of itself gives an equivalent algebra.

However, this explanation is not sufficient to explain why the algebra $T^*(\sl(2,\C))$ has a deformation along itself, because the cochain representing the algebra is trivial.  In this 
case, a nontrivial cocycle not given by the algebra determines a deformation along the algebra. 


\bibliographystyle{amsplain}

\providecommand{\bysame}{\leavevmode\hbox to3em{\hrulefill}\thinspace}
\providecommand{\MR}{\relax\ifhmode\unskip\space\fi MR }
\providecommand{\MRhref}[2]{%
  \href{http://www.ams.org/mathscinet-getitem?mr=#1}{#2}
}
\providecommand{\href}[2]{#2}

\end{document}